\documentclass{amsart}
\usepackage{amscd,amssymb}

\theoremstyle{plain}
\newtheorem{thm}{Theorem}[section]
\newtheorem{lem}[thm]{Lemma}
\newtheorem{definicja}[thm]{Definition}

\newtheorem{prop}[thm]{Proposition}

\theoremstyle{definition}

\theoremstyle{remark}

\newtheorem{Assumption}[thm]{Assumption}

\numberwithin{equation}{section}

\begin{document}
\setcounter{section}{-1}

\title{Non-nesting actions of Polish groups on pretrees}
\author{A.Ivanov}
\footnote{The research is supported by KBN grant 2 P03A 007 19}
\footnote{{\em E-mail address}: ivanov@math.uni.wroc.pl}
\footnote{{\em Fax number}: 48-71-3757429}

\maketitle
\centerline{Institute of Mathematics, Wroc{\l}aw University, pl.Grunwaldzki 2/4, }
\centerline{50-384, Wroc{\l}aw, Poland}
\bigskip
{\bf Abstract}
\begin{quote}
We study non-nesting actions of groups on
$\mathbb{R}$-trees.
We prove some fixed point theorems for such actions
under the assumption that a group is Polish and has
a comeagre conjugacy class.\parskip0pt


\end{quote}

\section{Introduction}

Non-nesting actions by homeomorphisms on $\mathbb{R}$-trees
frequently arise in geometric group theory (when actions on
spaces more general than trees are considered).
Explicitly they were introduced in \cite{levitt}.
We concentrate on the question when a Polish group
$G$ has fixed points under non-nesting
actions on $\mathbb{R}$-trees.
In Section 3 we prove the following theorem:
\begin{quote}
Let a Polish group $G$ have a non-nesting action on an
$\mathbb{R}$-tree $T_{0}$ without $G$-fixed points in $T_{0}$.
Let $X \subseteq G$ be a comeagre set.
Then the following statements hold. \\
If every element of $X$ fixes a point, then every element of $G$
fixes a point. \\
If $G$ fixes an end and $X$ is a conjugacy class, then
every element of $G$ fixes a point.
\end{quote}

These results are related to the paper of D.Macpherson 
and S.Thomas \cite{mactho} where they study actions of 
Polish groups on simplicial trees. 
Moreover in Section 3 we generalize the main result of 
\cite{mactho} that if a Polish group has a comeagre conjugacy 
class then every element of the group fixes a point under 
any action on a $\mathbb{Z}$-tree without inversions.
Our generalization concerns a wide class of non-nesting
actions on pretrees covering the case of isometric actions 
on $\mathbb{R}$-trees.
We apply this to actions of the group $Sym(\omega)$.\parskip0pt

On the other hand we are not able to extend the theorem of 
Macpherson and Thomas to non-nesting actions on 
$\mathbb{R}$-trees in general.
The main difficulty is that in this case we lose 
several basic properties of isometric actions (for example 
based on the material of Section 1 of \cite{cullmorg}).  
To remedy the situation we have made some general 
investigation of non-nesting actions on $\mathbb{R}$-trees.
This material is contained in Sections 1 and 2. \parskip0pt

In these sections we apply some axiomatic approach and as 
a result we really study a more general class of actions.  
First of all it turns out that the most appropriate
language for actions on trees by homeomorphisms is that 
of the betweenness relation $B(x;y,z)$; the corresponding 
structures are called {\em pretrees}.
The author has introduced in \cite{ivanov}
{\em classical actions on pretrees} and has noticed there that
$\mathbb{R}$-trees with isometric $G$-actions are classical.
Moreover in Section 2.2 below we show that non-nesting actions
by homeomorphisms on $\mathbb{R}$-trees are classical too.
The most interesting thing is that generalizing the theorem 
of Macpherson and Thomas we use methods which work for 
classical actions on median pretrees in general.
In particular we apply some new statement about
products of loxodromic elements (Proposition \ref{3loxo} below)
which can be considered as a metric-free version of 
sections 1.5 - 1.11 from \cite{cullmorg}.
\parskip0pt 

It is also worth noting that in our considerations we really 
use some algebraic property of comeagre classes 
(Condition (1) of Proposition \ref{4.1});
thus the results can be formulated in elementary terms.
\parskip0pt

Trying to generalize the theorem of Macpherson and Thomas to 
non-nesting actions on $\mathbb{R}$-trees we cannot eliminate 
the case when the comeagre conjugacy class consists of loxodromic 
elements. 
This case is investigated in Section 4 where we prove the most 
complicated result of the paper.
It roughly says that the presence of a comeagre conjugacy class
of loxodromic elements implies that orbits of any end-stabilizer
are much smaller that the corresponding orbits of the group.
It is also based on Proposition \ref{3loxo}.\parskip0pt 

Several preliminary versions of this paper have appeared since 1999.
A very close material is contained in Section 8 of the recent 
paper \cite{rosendal} of Rosendal, were the theorem of 
Macpherson and Thomas is extended to $\Lambda$-trees.        
Papers \cite{kecros} and \cite{solros} also study 
Polish groups with comeagre conjugacy classes and their 
actions on metric spaces.   
It is worth noting that some related problems 
(for example of embedding of generalized trees into 
$\mathbb{R}$-trees and $\Lambda$-trees) have been studied before 
(see \cite{bow}, \cite{chi}, \cite{dun} and \cite{levitt}).
Our motivation is partially based on these investigations.

\section{Median pretrees}

In this section we develop our basic tools which we will
later apply to the main theorems of the paper.
We start with very general results held
for group actions on median pretrees.\parskip0pt

The following definitions are taken from \cite{ivanov}.
Basically they appear in \cite{bow}.
The definition of a pretree is related to the definition
of a B-relation given in \cite{an}.

\begin{definicja}
A ternary structure $(T,B)$ is a pretree if the following axioms
are satisfied: \parskip0pt

* $(\forall x,y,z)(\neg B(y;x,x) \wedge \neg(B(y;x,z) \wedge B(z;x,y)))$;

* $(\forall x,y,z)(B(y;x,z) \leftrightarrow B(y;z,x))$;

* $(\forall x,y,z,w)(B(z;x,y) \wedge z \not= w \rightarrow
(B(z;x,w) \vee B(z;y,w)))$;
\end{definicja}

Define $[t,t'] = \{ x \in T: B(x;t,t') \vee x = t \vee x = t'\}$
the (closed) {\em interval (segment)} with endpoints $t,t'$.
We say that $[t,t')$ (and $(t,t']$, $(t,t')$ under the natural
definition) is an interval too.
A nonempty subset $S \subseteq T$ is an {\em arc}, if $S$ is
{\em full} (that is $(\forall x,y \in S)[x,y] \subseteq S$) and
{\em linear} (for all distinct $x,y,z \in S$ we have
$B(y;x,z) \vee B(z;x,y) \vee B(x;y,z)$). \parskip0pt

A pretree is {\em complete} if every arc is an interval,
not necessarily closed.

A point $x \in T$ is {\em terminal}, if
$(\forall y,z \in T)\neg B(x;y,z)$.
The pretree $T$ can be naturally decomposed
$T = T_{0} \cup P$, where $P$
is the set of all terminal points.  \parskip0pt

The pretree $(T,B)$ is {\em median} if for any $x,y,z \in T$
there is an element  $c \in [x,y]\cap [y,z] \cap [z,x]$.
In this case $c$ is unique and is called the {\em median}
of $x,y,z$; we will write $c = m(x,y,z)$. \parskip0pt

\begin{Assumption}
{\em From now on we consider only median pretrees.}
\end{Assumption}

It is clear that every simplicial or real tree can be considered
as a complete median pretree by adding ends as terminal points and 
taking the reduct to the natural betweenness relation.
\parskip0pt

The following notion will be applied below several times (for
subpretrees of median pretrees).
It has not been formulated before.
We say that a pretree is {\em quasimedian} if for any triple
$t,q,r$, if the interval $[t,q)$ is not closed and is contained
in $[t,r]$ then $q\in [t,r]$.
To see that {\em a median pretree is quasimedian} let $c=m(t,q,r)$.
If $c=q$ then $q\in [t,r]$.
If $c\not= q$, then since $[t,q)$ is not closed and $c\in [t,q)$,
the interval $(c,q)=[t,q)\setminus [t,r]$ is not empty,
contradicting the assumptions.
\parskip0pt

Let $x \in T_{0}$.
A maximal arc of the form $L_{x} =\bigcup [x,t_{\gamma}]$
where all $t_{\gamma}$ are not terminal, is called a {\em half-line}.
Half-lines $L_{x}$ and $L_{y}$ are equivalent
if there is a half-line $L_{z} \subseteq L_{x} \cap L_{y}$.
An {\em end} ${\bf e}$ of $T$ is an equivalence class of
half-lines.
Define a partial order $<_{\bf e}$ by $x <_{\bf e} y$
if the half-line $L_{x}$ representing ${\bf e}$ contains $y$.
\parskip0pt

It is clear that an arc of the form $[x,p), p \in P,$
is a half-line.
Since $p$ is a terminal point, any pair of half-lines $[x,p)$
and $[y,p)$ with $x,y\in T_0$, have a common point from $T_{0}$
(which is the corresponding median).
This shows that the set of all half-lines $[t,p), t\in T_{0}$
forms an end (the {\em end corresponding} to $p\in P$).\parskip0pt

A maximal arc of the form $\bigcup [t_{\gamma},t'_{\gamma}]$
where $t_{\gamma},t'_{\gamma}$ are not terminal, is called a
{\em line}.
It is worth noting that the ends of a line of a complete pretree $T$ 
are presented by a pair of termial points of $T$.
The following lemma is conceivably known (see \cite{tits},\cite{bc}
and Section 2 of \cite{bow}) and is based on existence of medians.
A complete proof of the lemma (in a slightly more general situation)
is given in \cite{ivanov}.

\begin{lem} \label{4sts}
1. The intersection of a finitely many (distinct) segments or
half-lines with a common extremity $t \in T_{0}$ is a segment
having $t$ as an extremity. \parskip0pt

2. If $t,q,r \in T_{0}$ satisfy $[t,q] \cap [q,r] = \{ q\}$, then
$[t,r] = [t,q] \cup [q,r]$.
If the interval $[t,q)$ is not closed and $[t,q)=[t,r)$, then
$q = r$. \parskip0pt

3. Let $T$ be complete and $A$ and $B$ be two full subsets of $T_{0}$ 
whose intersection consists of at most one point.
Then there exists a segment $[t,q]$ which is contained in every
full set which has a non-empty intersection with both $A$ and $B$,
and, moreover, $A\cup \{ t\}$ and $B\cup \{ q \}$ are full.
If for some $\epsilon,\tau \in \{ 0,1\}$, $t$ and $q$ satisfy
$t\in^{\epsilon} A \wedge q\in^{\tau} B$ (where $\in^{0}$ denotes
$\not\in$), then the condition
$t\in^{\epsilon} A\wedge q\in^{\tau} B$ determines the segment
uniquely. \parskip0pt

4. Under the circumstances of the previous statement if $A$ is a segment, 
a line or a half-line, then $t$ as above can be found in $A$.
Moreover in the case when $A$ and $B$ are lines these statements hold 
without the assumtion that $T$ is complete. 
\end{lem}

Let $T$ and $A,B\subseteq T_{0}$ with $A\cap B=\emptyset$, satisfy 
the assumptions of Lemma \ref{4sts} (3 or 4). 
Then there exists a segment $[t,q]$ which is
contained in every full set which has a non-empty intersection
with both $A$ and $B$, and, moreover, $A\cup\{ t\}$ and
$B\cup\{ q\}$ are full.
Then we call the interval $[t,q]\setminus (A\cup B)$ the
{\em bridge between} $A$ and $B$.
It can happen that the bridge is an open (or empty) interval.
In this case we define the bridge by its extremities $t$ and $q$ 
as $(t,q)$.
\parskip0pt

Let $G$ be a group acting on a median pretree
$T$ by automorphisms of the structure $(T,B)$.
It is clear that the set $P$ of terminal points is $G$-invariant.
Let $g \in G$.
The set of $g$-fixed points is denoted by $T^{g}$.
The element $g\in G$ is {\em loxodromic}, if
$T^{g}_{0} =\emptyset$, $|g| = \infty$ and there exists a unique
$g$-invariant line in $T_{0}$ such that $g$ preserves
the natural orders on the line.
It is called the {\em axis (characteristic line)} of $g$.
In the case of an isometric action on an $\mathbb{R}$-tree
a loxodromic element is {\em hyperbolic}.
The proof of the following lemma is standard
(by arguments from \cite{tits}, Section 3.1)
and can be found in \cite{ivanov}.

\begin{lem} \label{loxo}
Let $G$ act on $T$ and $g \in G$.
If $g$ is loxodromic, then \parskip0pt

(a) for any $p \in T_{0}$ the segment
$[p,g(p)]$ meets the characteristic line $L_{g}$
and $[p,g(p)] \cap L_{g} = [q,g(q)]$ for some $q \in L_{g}$.
\parskip0pt

(b) $x \in L_{g}$ if and only if $x$ is the median of
$x,g^{-1}(x),g(x)$.
\end{lem}

The following proposition can be considered as a metric-free 
version of statements 1.5 - 1.11 of \cite{cullmorg}.
It is very convenient in applications and
will be one of the basic tools of the paper.

\begin{prop} \label{3loxo}
Let $G$ act on a median pretree $T$.
Let $h_{1},h_{2},h_{3} \in G$ be loxodromic and
$h_{2} \cdot h_{1} = h_{3}$.
Then one of the following cases holds:\parskip0pt 

(1) there are $d \in L_{h_{1}}$, $b \in L_{h_{2}}$ and
$c \in L_{h_{3}}$ such that one of the segments
$[d,h_{1}(d)]$, $[b,h_{2}(b)]$ and $[c,h_{3}(c)]$ properly
contains the others;\parskip0pt

(2) there is a non-linear triple $d,b,c$ with 
$m(d,b,c)\in L_{h_1}\cap L_{h_2}\cap L_{h_3}$ such that 
$c\in L_{h_2}\cap L_{h_3}$, $d\in L_{h_1}\cap L_{h_3}$,
$[m(d,b,c),b]= L_{h_1}\cap L_{h_2}$ and  
$h_1 (d)=b$, $h_2 (b)=c$, $h_3(d)=c$.
In this case the element $h_4 = h^{-1}_2 h_1$ has no 
fixed points and the set
$\bigcup \{ [h^{j}_4 (d),h^{j+1}_4 (d)]:j\in \mathbb{Z}\}$ 
is an arc such that its segment $[d,h_4 (d)]$  
properly contains $[d,h_1 (d)]$ and $[b,h^{-1}_2 (b)]$. 
\end{prop}

{\em Proof.}
Consider the case when $|L_{h_{1}}\cap L_{h_{2}}|\le 1$.
Let a segment $[a,b]$ satisfy Lemma \ref{4sts}(3,4) with 
respect to $L_{h_{1}}$ and $L_{h_{2}}$ (thus $a \in L_{h_{1}}$ 
and $b\in L_{h_{2}}$).
Then applying Lemma \ref{4sts}(2) three times,
$$
[a,h_{2}(h_{1}(a))]=[a,b]\cup [b,h_{2}(b)]\cup [h_{2}(b),h_{2}(a)]
\cup [h_{2}(a),h_{2}(h_{1}(a))].
$$
Moreover, by Lemma \ref{loxo} the segment $[a,h_{2}(h_{1}(a))]$
meets $L_{h_{3}}$.
Similarly the segment $[h^{-1}_{1}(b),h_{2}(b)]$ meets
$L_{h_{3}}$ and
$$
[h^{-1}_{1}(b),h_{2}(b)] = [h^{-1}_{1}(b),h^{-1}_{1}(a)] \cup
[h^{-1}_{1}(a),a] \cup [a,b] \cup [b,h_{2}(b)].
$$
Since $[h^{-1}_{1}(b),h_{2}(h_{1}(a))] =$
$$
[h^{-1}_{1}(b),h^{-1}_{1}(a)]\cup [h^{-1}_{1}(a),a]\cup [a,b]
\cup [b,h_{2}(b)]\cup [h_{2}(b),h_{2}(a)]\cup
[h_{2}(a),h_{2}(h_{1}(a))],
$$
the line $L_{h_{3}}$ must meet $[a,b]\cup [b,h_{2}(b)]$.\parskip0pt

If $[b,h_{2}(b)] \cap L_{h_{3}} \not= \emptyset$, then
$[h^{-1}_{1}(h^{-1}_{2}(b)),h^{-1}_{1}(b)]\cap L_{h_{3}}\not=\emptyset$
and the elements between those intersections lie in $L_{h_{3}}$.
Thus $h^{-1}_{1}(b),h^{-1}_{1}(a),a,b \in L_{h_{3}}$.
This implies
$h_{2}(h_{1}(a)) \in L_{h_{3}}$ and $h_{2}(b) \in L_{h_{3}}$.
\parskip0pt

If $[a,b] \cap L_{h_{3}} \not= \emptyset$, then
$[h_{2}(h_{1}(a)),h_{2}(h_{1}(b))]\cap L_{h_{3}}\not= \emptyset$
and the elements between those intersections lie in $L_{h_{3}}$.
Thus $b,h_{2}(b)\in L_{h_{3}}$.
This implies $h^{-1}_{1}(b)\in L_{h_{3}}$
and $h^{-1}_{1}(a), a \in L_{h_{3}}$. \parskip0pt

Let $d := h^{-1}_{1}(a)$ and $c:= h^{-1}_{1}(b)$.
Then $[b,h_{2}(b)]$ and $[d,h_{1}(d)]$ are proper subsegments
of $[c,h_{3}(c)] = [h^{-1}_{1}(b),h_{2}(b))]$. \parskip0pt

If $|L_{h_{1}} \cap L_{h_{2}}| \ge 1$, let $e_1 ,e'_{1}$ be 
the ends of (half-lines of) $L_{h_{1}}$ and $e_2 ,e'_{2}$ be 
the ends of $L_{h_{2}}$. 
For ease of notation we extend the betweenness relation of 
$T_0$ to $T_0 \cup \{ e_1 ,e'_{1},e_2 ,e'_{2} \}$ in the 
obvious way (so that $L_{h_{1}} = [e_{1},e'_{1}] \cap T_{0}$ 
and $L_{h_{2}} = [e_{2},e'_{2}] \cap T_{0}$).
\parskip0pt 

If the lines $L_{h_{1}}$ and $L_{h_{2}}$ represent the same 
end $e \in \{ e_{1},e'_{1}\}\cap\{ e_{2},e'_{2}\}$,
then $L_{h_{3}}$ represents $e$.
Thus there exists
$d \in L_{h_{1}} \cap L_{h_{2}} \cap L_{h_{3}}$ with
$$
h_{1}(d),h^{-1}_{1}(d),h_{2}(d),h^{-1}_{2}(d),h_{3}(d),
h^{-1}_{3}(d)\in L_{h_{1}} \cap L_{h_{2}} \cap L_{h_{3}}.
$$
Now the lemma is obvious for $b := h_{1}(d)$
and $c = d$. \parskip0pt

If $\{ e_{1},e'_{1}\}\cap\{ e_{2},e'_{2}\} = \emptyset$,
find $a',a'' \in L_{h_{1}} \cap L_{h_{2}}$ such that
$[a',a''] = L_{h_{1}}\cap L_{h_{2}}$ (notice that the points
$a'$ and $a''$ can be found as medians of appropriate triples
from $\{ e_{1},e'_{1},e_{2},e'_{2}\}$).
Assume that there are $q_{1}\in\{ e_{1},e'_{1}\}$
and $q_{2}\in\{ e_{2},e'_{2}\}$, such that
$a' = m(q_{2},a'',q_{1})$ and the following condition holds:
$(h^{-1}_{1}(a')\in [a',q_{1}])$ or $(h_{2}(a')\in [a',q_{2}])$.
Then $a'\in [h_{2}(a'),h^{-1}_{1}(a')]$.
When the conjunction 
$(h^{-1}_{1}(a')\in [a',q_{1}])\wedge (h_{2}(a')\in [a',q_{2}])$
does not hold we put $b = a := a'$ and apply the argument of 
the case $|L_{h_{1}} \cap L_{h_{2}}| \le 1$ 
(with some simplifications) to an appropriate pair of the segments
$[h^{-1}_{1}(h^{-1}_{2}(a)),a]$, $[h^{-1}_{1}(b),h_{2}(b)]$ and
$[h^{-1}_{1}(b),h_{2}(b)]$, $[a,h_{2}(h_{1}(a))]$. 
This argument shows that the segments 
$[h^{-1}_{1}(a'),a']$ and $[a',h_{2}(a')]$ lie on $L_{h_3}$ as 
proper subsegments of $[c,h_{3}(c)]$ where $c=h^{-1}_{1}(a')$.
\parskip0pt 

The case of the conjunction 
$(h^{-1}_{1}(a')\in [a',q_{1}])\wedge (h_{2}(a')\in [a',q_{2}])$
is divided into two subcases. 
The first one appears if $h_1 (a')\in [a',a'']$ or 
$h^{-1}_2 (a')\in [a',a'']$. 
Then applying arguments as above we obtain that $a'\in L_{h_3}$ and 
either 
$[a',h_3 (a')]\cup [a',h_1 (a')]\subseteq [h_1 (a'), h_2 (h_1 (a'))]$
or 
$[h^{-1}_1 (h^{-1}_2 (a'),a']\cup [h^{-1}_2 (a'),a']\subseteq [h^{-1}_1 (h^{-1}_2 (a'),h^{-1}_2 (a')]$. 
\parskip0pt

In the case of the conjunctions 
$(h^{-1}_{1}(a')\in [a',q_{1}])\wedge (h_{2}(a')\in [a',q_{2}])$
and $h_1 (a')\not\in [a',a''] \wedge h^{-1}_2 (a')\not\in [a',a'']$
let $b=a''$, $d=h^{-1}_1 (a'')$ and $c=h_2 (a'')$.
It is easy to see that the second condition of the statement 
of the proposition is satisfied (where $a'=m(d,b,c)$). 
It is straightforward that for all $j\in \mathbb{Z}$,  
$[h^{j-1}_4 (d),h^{j}_4 (d)]\cap [h^{j}_4 (d),h^{j+1}_4 (d)]=\{ h^{j}_4 (d)\}$; 
thus $\bigcup \{ [h^{j}_4 (d),h^{j+1}_4 (d)]:j\in \mathbb{Z}\}$ is an arc. 
To see that $h_4$ does not fix any point assume that $h^{-1}_2 (h_1 (v))=v$.
Then the segment $[v,h_1 (v)]$ meets both $L_{h_1}$ and 
$L_{h_2}$ and thus $[a',a'']$. 
Let $[v,h_1 (v)]\cap L_{h_1}=[w,h_1 (w)]$. 
The assumption  $h_1 (a')\not\in [a',a'']$ implies that 
one of the extremities $w$ or $h_1 (w)$ is outside of $[a',a'']$. 
If $w\not\in [a',a'']$ the assumtion $h^{-1}_2 (a')\not\in [a',a'']$ 
now implies that $[w,a']\cup [a',a'']$ is a non-trivial 
subsegment of $[v,h^{-1}_2 (h_1 (v))]$, a contradiction
with $h^{-1}_2 (h_1 (v))=v$. 
If $h_1 (w)\not\in [a',a'']$ then similarly 
$[a'',h^{-1}_2 (a'')]\cup [h^{-1}_2 (a''),h^{-1}_2 (h_1 (w))]$ 
is a non-trivial subsegment of $[v, h^{-1}_2 (h_1 (v))]$, 
a contradiction.  \parskip0pt 

Now assume that
$\{ e_{1},e'_{1}\}\cap\{ e_{2},e'_{2}\} = \emptyset$,
but none of the cases 
$(h^{-1}_{1}(a')\in [a',q_{1}])$ and $(h_{2}(a')\in [a',q_{2}])$
holds where $[a',a''] = L_{h_{1}}\cap L_{h_{2}}$ 
and $q_1$ and $q_2$ are as above.
Now we have that $h_{1}(a') \in [a',q_{1}]$ and 
$h^{-1}_{2}(a')\in [a',q_{2}]$.
We may assume that $q_{1} = e_{1}$ and $q_{2} = e_{2}$.
Since $\{ e_{1},e'_{1}\}\cap\{ e_{2},e'_{2}\} = \emptyset$,
we have $a''= m(e'_{2},a',e'_{1})$, $h_{1}(a'')\in [a'',e'_{1}]$
and $h^{-1}_{2}(a'')\in [a'',e'_{2}]$ (otherwise the arguments of 
the paragraphs above work for $a''$ instead of $a'$).
Then $h^{-1}_{1}(a'),h_{2}(a')\in [a',a'']$ and
$h^{-1}_{1}(a''),h_{2}(a'')\in [a',a'']$ (because 
$h_{1}$ and $h_{2}$ preserve the natural orderings of their axises).
If $h^{-1}_{1}(a')\in L_{h_{3}}$ then taking $d = c = h^{-1}_{1}(a')$ 
and $b = a'$ we see that one of the segments  $[d,h_{1}(d)]$ or 
$[b,h_{2}(b)]$ properly contains the remaining ones 
(including $[c,h_{3}(c)]$).
The same argument works or the case $h^{-1}_{1}(a'')\in L_{h_{3}}$.
\parskip0pt

Let $c' = h^{-1}_{1}(a')$ and $c'' = h^{-1}_{1}(a'')$.
We know that
$[c',h_{3}(c')]\cup [c'',h_{3}(c'')]\subseteq [a',a'']$.
Since $[c',h_{3}(c')]$ and $[c'',h_{3}(c'')]$ meet $L_{h_{3}}$,
in the case $[c',h_{3}(c')] \cap [c'',h_{3}(c'')] =\emptyset$
we have $c',c'' \in L_{h_{3}}$ and by the previous paragraph,
this is enough for the proposition.
Assume $[c',h_{3}(c')] \cap [c'',h_{3}(c'')] \not=\emptyset$.
Since $h_{1}$ and $h_{2}$ preserve the natural orders of their 
axises, $c'\in [c',h_{3}(c')] \cap [c'',h_{3}(c'')]$ or
$c'' \in [c',h_{3}(c')] \cap [c'',h_{3}(c'')]$.
We may assume that $c'\in [c',h_{3}(c')] \cap [c'',h_{3}(c'')]$
(then $h_{3}(c'') \in [c',h_{3}(c')]$).
Since $[c',h_{3}(c')]$ and $[c'',h_{3}(c'')]$ meet $L_{h_{3}}$,
there is $u \in L_{h_{3}} \cap [c', h_{3}(c'')]$.
Then $h^{-1}_{2}(u) \in [e'_{2},a'']$ and
$[h^{-1}_{1}(h^{-1}_{2}(u)),c''] \cap [c'',a'] = \{ c''\}$
(because
$[h^{-1}_{1}(h^{-1}_{2}(u)),c'']=h^{-1}_{1}([h^{-1}_{2}(u)),a''])$).
This implies $c',c''\in [h^{-1}_{3}(u),u]$.
By $h^{-1}_{3}(u),u\in L_{h_{3}}$ we have $c',c''\in L_{h_{3}}$,
which is enough for the proposition.
$\square$

\section{Classical actions on pretrees}

This section contains preliminary results on non-nesting
classical actions.

\subsection{Classical actions and non-nesting actions}

In the following definition we collect usual properties of
typical actions  (for example {\em isometric} ones).

\begin{definicja} \label{defcl}
The action of $G$ on a complete median pretree $T$ 
($=T_{0} \cup P$) is {\em classical} if
the following conditions hold: \parskip0pt

(C0)  If $x,y \in T^{g}_{0}$ ($:= T^{g} \cap T_{0}$), then
$[x,y] \subseteq T^{g}_{0}$; \parskip0pt

(C1) If $x \not\in T^{g}_{0} \not= \emptyset$, then
$|[x,g(x)] \cap T^{g}_{0}| = 1$; \parskip0pt

(C2)  If $T^{g}_{0} = \emptyset$, then $g$ is
loxodromic;\parskip0pt

(C3)  If $g_{1},g_{2},g_{3} \in G$ are not loxodromic and
$g_{1}=g_{2}\cdot g_{3}$, then
$T^{g_{1}}_{0} \cap T^{g_{2}}_{0} \cap T^{g_{3}}_{0} \not= \emptyset$; \parskip0pt

(C4)  If $g$ is loxodromic, $L$ is the axis of $g$ and
$g = h\cdot h'$ with $T^{h}_{0}\not= \emptyset\not= T^{h'}_{0}$, then
$|T^{h}_{0} \cap L| = 1 = |T^{h'}_{0}\cap L|$.
\end{definicja}

It is worth noting here that we do not really need condition
(C1) in this paper.
We include it into the definition of classical actions because
non-nesting actions satisfy it (see below).
On the other hand any pretree satisfies a weaker form of
$(C1)$: if $x \not\in T^{g}_{0} \not= \emptyset$, then
$|[x,g(x)] \cap T^{g}_{0}| \le 1$. Indeed,
if $u,u' \in T^{g}_{0}\cap [x,g(x)]$ and $u \in [x,u')$,
then $[x,g(x)] = [x,u') \cup [u',g(x)]$ and
$u = g(u)\in g([x,u']) = [g(x),u']$, a contradiction. \parskip0pt

Let $\Lambda$ be linearly ordered abelian group.
A $\Lambda$-metric space $(X,d)$ is called
a $\Lambda$-{\em tree} if:\parskip0pt

(1) $X$ is {\em geodesically linear}: for any $a,b \in X$ there
exists a unique metric morphism $\alpha :[0,d(a,b)]\rightarrow X$
such that $\alpha(0) =a$ and $\alpha(d(a,b)) = b$
(then $[a,b] := \alpha([0,d(a,b)])$).\parskip0pt

(2) $\forall x,y,z\exists !w([x,y]\cap [x,z]=[x,w])$. \parskip0pt

(3) 
$\forall x,y,z([x,y]\cap [y,z]=\{ y\}\rightarrow [x,y]\cup [y,z] = [x,z])$.
\parskip0pt

For $\mathbb{R}$-trees and isometric actions on them
Lemmas \ref{4sts} and \ref{loxo} are known \cite{tits}
(moreover in Lemma \ref{4sts}(3) for closed $A$ and $B$ we have
$t \in A$ and $q\in B$).
It is also known that every action of a group on an
$\mathbb{R}$-tree by isometries induces a classical action on
the pretree extended by all ends with respect to the natural
betweenness relation.
The proof can be extracted from \cite{tits} (a more convenient
reference is \cite{cullevogt}, where Lemma 1.2 corresponds to
propery (C4) of Definition \ref{defcl}).
It is worth noting that $T^{g}_{0}$ is a closed set if $g$ is
an isometry of an $\mathbb{R}$-tree $T_{0}$. \parskip0pt

The following axiom ({\em non-nesting}) defines classical
actions quite close to isometric ones.
$$
\forall g\in G \forall t,t'\in T_{0} \neg ([g(t),g(t')]\subset
[t,t']).
$$

The axiom of non-nesting has the following immediate consequences.
(The following lemma appears in a different form in \cite{bc}.)

\begin{lem} \label{1.9}
Let a group $G$ have a classical non-nesting action on a pretree
$T$ and let $g \in G$ be loxodromic.
Then under an appropriate choice
of $+\infty$ on $L_{g}$, the element $g$ is strictly increasing
on $L_{g}$.
\end{lem}

{\em Proof.} Let $a,b \in L_{g}$ and $g(a) < a < b \le g(b)$.
Then $g^{-1}$ maps $[g(a),g(b)]$ properly into itself.
The remaining cases are similar.$\square$

\bigskip

Under assumptions of the lemma let $L$ be a line of $T$ and
$G_{\{ L\} }$ (and $G_{L}$) be the stabilizer
(pointwise stabilizer) of $L$ in $G$.
Assume that $G_{\{ L\} }$ does not have elements reversing
the terminal points of $L$
(this happens when it does not have subgroups of index 2).
Then non-nesting implies that given an ordering of $L$
corresponding to the betweenness relation of $T$, the following
relation makes $G_{\{ L\} }/G_{L}$ a linearly ordered group:
$g<g'$ if and only if $\exists t\in L (g(t)< g'(t))$. \parskip0pt

Now assume that $H$ is a subgroup of the stabilizer $G_{\{ L\} }$.
By completeness $L$ consists of $H$-invariant intervals.
If $[a,b)$ is such an interval then by non-nesting,
$h(a) \not\in [a,b)$, where $h\in H$.
We see that $H$ acts trivially on $L$ or
there are no proper $H$-invariant arcs in $L$.
\parskip0pt

This argument also shows that for a loxodromic $g$ the axis
$L_{g}$ does not have proper $g$-invariant subintervals.
In the case when $L_{g}$ is homeomorphic with $\mathbb{R}$ we
immediately have that
{\em up to topological conjugacy $g$ can be viewed
as a translation by a real number.}
On the other hand we also have that {\em for every
$h\in G_{\{ L_{g}\} }$ there exists $n\in \omega$ such that}
$h<g^{n}$.

\subsection{Non-nesting actions on $\mathbb{R}$-trees 
and their end stabilizers}

Let $G$ be an infinite group acting on an $\mathbb{R}$-tree
$T_{0}$ by homeomorphisms 
\footnote{preserving the betweenness relation}.
By $T$ we denote $T_{0}$ together with the set of ends.
As above the action is called {\em non-nesting}
\cite{levitt} if no $g \in G$ maps an arc properly into itself.
In this case if $T^{g}_{0}$ is not empty then $T^{g}_{0}$ is
closed.
If $g$ does not fix any point, then by Theorem 3(2) from
\cite{levitt} there exists a geodesic $\mathbb{R}$-line
$L_{g}\subseteq T_{0}$ such that $g$ acts on $L_{g}$ by 
an order preserving transformation, which is a translation 
up to topological conjugacy.
The following proposition develops these observations.

\begin{prop} \label{3.1}
Let $G$ be a group.
Every non-nesting action of $G$ on an $\mathbb{R}$-tree $T_{0}$
induces a classical action on the corresponding pretree.
\end{prop}

{\em Proof.}
Since an $\mathbb{R}$-tree is complete and median we must
verify (C0)-(C4) of Definition \ref{defcl}.
Condition (C0) is clear.
Conditions (C1) and (C2) are proved in Theorem 3(1,2) of
\cite{levitt}. \parskip0pt

To see (C3) let each of $h_{1},h_{2},h_{3} \in G$ fix points in
$T_{0}$ and $h_{2} \cdot h_{1} = h_{3}$.
We want to show that
$T^{h_{1}}_{0}\cap T^{h_{2}}_{0}\cap T^{h_{3}}_{0}\not=\emptyset$.
\parskip0pt

Take $t\in T^{h_{1}}_{0}$ with minimal distance
$d(t,T^{h_{2}}_{0})$ (with minimal distance from $T^{h_{2}}_{0}$).
Then the segment $[t,h_{2}(t)]$ meets $T^{h_{2}}_{0}$
at precisely one point, say $c$.
Since $h_{2}(t) =h_{3}(t)$, $[t,h_{2}(t)]$ meets $T^{h_{3}}_{0}$
at precisely one point too, say $d$.
If these points are not the same then one of the maps
$h^{-1}_{2}h_{3}$ or $h^{-1}_{3}h_{2}$ maps an arc $[t,d]$ or
$[t,c]$ into itself (the first case corresponds to the situation
when $c$ is between $t$ and  $d$).
So, there is a point fixed by $h_{2}$ and $h_{3}$.
It must be fixed by $h_{1}$.\parskip0pt

The proof of (C4) is related to Proposition 1.2 from
\cite{cullevogt}. \parskip0pt

Let $g = h'h$ be loxodromic, but $h,h' \in G$
fix some points in $T_{0}$.
Let $a \in T^{h}_{0}$, $b \in T^{h'}_{0}$ be chosen with minimal
$d(a,b)$ (the existence of such $a$ and $b$ follows from the
fact that $T^{h}_{0}$ and $T^{h'}_{0}$ are closed).
Let us prove that
$L_{g}\cap T^{h}_{0}=\{ a\}$ and $L_{g}\cap T^{h'}_{0}=\{ b\}$.
\parskip0pt

The segment $[a,g(a)] = [a,h'(a)]$ meets $L_{g}$ and contains
$b$ (by (C1) and Lemma \ref{loxo}).
Assume $[a,g(a)]= [a,q]\cup [q,g(q)]\cup [g(q),g(a)]$
where $[q,g(q)] = L_{g} \cap [a,g(a)]$.
If $b \not\in L_{g}$, for example
$b\in [a,q)$, then $g(h')^{-1}([b,g(a)]) = [g(b),g(a)]$ is
properly contained in $[b,g(a)]$ (because
$[b,g(a)]= [b,q]\cup [q,g(q)]\cup [g(q),g(b)]\cup [g(b),g(a)]$).
The case when $b$ is between $g(a)$ and an element
from $L_{g}$ is similar.
Thus $b \in L_{g}$ and $h^{-1}(b) = g^{-1}(b) \in L_{g}$.
Then the  segment $[g^{-1}(b),b]$ belongs to $L_{g}$.
Thus $a \in L_{g}$. If there exists
$a' \in T^{h}_{0} \setminus \{ a \}$ belonging to $L_{g}$,
then $[a',a] \subset L_{g} \cap  T^{h}_{0}$.
Thus $[g^{-1}(b),b] \cap [a',a] = \{ a \}$ and we see that
$a',b,g^{-1}(b)$ belong to $L_{g}$ but are not linear, 
a contradiction. \parskip0pt

The proof that $L_{g}\cap T^{h'}_{0} =\{ b\}$ is similar.
$\square$
\bigskip

In a sense Lemma \ref{1.9} describes elements stabilizing
a line under the assumption of non-nesting.
We now concentrate on {\em end stabilizers of non-nesting
actions on $\mathbb{R}$-trees}.
We will see that some kind of Lemma \ref{1.9} holds in this case.
Let a group $G$ have a non-nesting action on
an $\mathbb{R}$-tree $T$ without $G$-fixed points in $T_{0}$.
Assume that there is a loxodromic $g \in G$.
Let $a_{0} \in L_{g}$ and ${\bf e}$ be the end represented by
$(-\infty ,a_{0}]$ ($-\infty$ is chosen so that $g$ is increasing).
Consider the stabilizer $G_{{\bf e}}$.
Let $G_{({\bf e})}$ be the subset of $G_{{\bf e}}$ of all
elements fixing some points in $T_0$. 
Note that any $h \in G_{{\bf e}}$ defines a map
$(-\infty ,a] \rightarrow (-\infty ,b]$ for some $a,b \le a_{0}$.
Thus $G_{({\bf e})}$ be the subgroup of elements fixing pointwise
$(-\infty ,a]$ for some $a \le a_{0}$ (by non-nesting).
We also see that it is normal in $G_{{\bf e}}$.
In the following lemma we consider the group
$G_{{\bf e}}/G_{({\bf e})}$.

\begin{lem} \label{1.10}
Let a group $G$ have a non-nesting action on an
$\mathbb{R}$-tree $T$ without $G$-fixed points in $T_{0}$.
Let $g \in G$ be loxodromic and ${\bf e}$ be the
$(-\infty)$-end of $L_{g}$.
Then the group $G_{{\bf e}}/G_{({\bf e})}$ is embeddable into
$(\mathbb{R},+)$ as a linearly ordered group under the ordering: 
$gG_{({\bf e})}\prec g'G_{({\bf e})}\Leftrightarrow \exists t\in T_0 (g'(t)<_{{\bf e}}g(t))$.
\end{lem}

{\em Proof.} By Lemma \ref{1.9} any
$h \in G_{{\bf e}}\setminus G_{({\bf e})}$ defines a map
$(-\infty ,a] \rightarrow (-\infty ,b]$ with some 
$a,b\in (-\infty ,a_0 ]$, which is strictly monotonic on $(-\infty ,a]$.
Now notice that under the induced ordering $\prec$ the group
$G_{{\bf e}}/G_{({\bf e})}$ is a linearly ordered group.
Indeed, linearity follows from Lemma \ref{1.9}.
If $g_{1},g_{2} \in G_{{\bf e}}$ satisfy $g_{1}(a) < g_{2}(a)$
for some $a \le a_{0}$, then by non-nesting
for all $a' \le a$, $g_{1}(a') < g_{2}(a')$.
We see that for every $g' \in G_{{\bf e}}$ there exists $b \le a$
such that $g_{1}\cdot g'(x) \le g_{2}\cdot g'(x)$ for all
$x \le b$.
On the other hand if $g_{1},g_{2} \in G_{{\bf e}}$ satisfy
$g_{1}G_{({\bf e})} \preceq g_{2}G_{({\bf e})}$ then obviously
$g'\cdot g_1 G_{({\bf e})}\preceq g'\cdot g_{2} G_{({\bf e})}$ 
for all $g' \in G_{{\bf e}}$.
This shows that $G_{{\bf e}}/G_{({\bf e})}$ is a linearly
ordered group. \parskip0pt

Since the elements of $G_{{\bf e}}$ act by translations up to
topological conjugacy, $G_{{\bf e}}/G_{({\bf e})}$ is Archimedean.
By H\"{o}lder's theorem it is a subgroup of $(\mathbb{R},+)$
\cite{birkhoff}.$\square$

\bigskip

{\em We now define an action $*_{g}$ of} $G_{{\bf e}}$
on $L_{g}$.
Let $h\in G_{{\bf e}}$ and $c \in L_{g}$.
Find a natural number $n_{0}$ such that $g^{-n_{0}}(c)$ is 
greater with respect to $<_{{\bf e}}$ than any of 
the elements $a_0 ,h(a_0 ), h^{-1}(a_0 )$.
Now let $h*_{g}c = g^{n_{0}}hg^{-n_{0}}(c)$.
By the choice of $n_{0}$ we see 
$h*_g c \in L_{g}$. \parskip0pt

It is worth noting that for every $n \ge n_{0}$,
$h*_{g}c = g^{n}hg^{-n}(c)$.
This follows from the fact
that the element $h^{-1}g^{n-n_{0}}hg^{n_{0}-n}$ belongs to
$G_{({\bf e})}$ (as $G_{{\bf e}}/G_{({\bf e})}$ is a subgroup of
$(\mathbb{R},+)$) and then (by non-nesting) the transformations
$hg^{n_{0}-n}$ and $g^{n_{0}-n}h$ are equal at $g^{-n_{0}}(c)$.
We now see:
$$
g^{n}hg^{-n}(c) =
g^{n_{0}}h(h^{-1}g^{n-n_{0}}\cdot hg^{n_{0}-n}(g^{-n_{0}}(c)))
=g^{n_{0}}hg^{-n_{0}}(c).
$$
Now it is easy to see that $*_{g}$ is an action and the elements
of $G_{({\bf e})}$ act on $L_{g}$ trivially. \parskip0pt

Let $L^{a_{0}}_{g}= G_{{\bf e}}a_{0} \cap L_{g}$ where
$G_{{\bf e}}a_{0}$ is the orbit of $a_{0}$.
Then there exists a surjection
$\nu_{a_{0}} : G_{{\bf e}}\rightarrow L^{a_{0}}_{g}$ defined
by $\nu_{a_{0}} (h) = h*_{g}a_{0}$ (with respect to the
action defined above).
It is easy to see that the map $\nu_{a_{0}}$ is surjective.
Moreover, for any $h,h'\in G_{{\bf e}}$,
$\nu_{a_{0}} (h\cdot h') = h*_{g}\nu_{a_{0}} (h')$.

\begin{lem}  \label{1.11}
The map $\nu_{a_{0}}$ defines an order-preserving bijection
from $G_{{\bf e}}/G_{({\bf e})}$ onto $L^{a_{0}}_{g}$ under
the order induced by $L_{g}$.
\end{lem}

{\em Proof.}
Notice that if $\nu_{a_{0}} (h_{1}) = \nu_{a_{0}} (h_{2})$
then $h_{1}h^{-1}_{2}$ fixes some $(-\infty ,a]$,
$a \le a_{0}$, pointwise.
Indeed, let $n$ and $a = g^{-n}(a_{0})$ be chosen so
that $a,h_{1}(a),h_{2}(a)$, $h^{-1}_{1}(a),h^{-1}_{2}(a)$
$\in L_{g} \cap L_{h_{1}} \cap L_{h_{2}}$
(for $h_i \in G_{({\bf e})}$ we replace $L_{h_i}$ by $T^{h_i}_0$).
Then
$$
g^{n}h_{1}g^{-n}(a_{0}) = \nu_{a_{0}} (h_{1}) =
\nu_{a_{0}} (h_{2}) = g^{n}h_{2}g^{-n}(a_{0})
$$
and we see that $h_{1}(a) = h_{2}(a)$.
Now the claim follows from non-nesting. \parskip0pt

Applying non-nesting again we obtain that the preimage of 
$a_{0}$ (with respect to $\nu_{a_{0}}$) equals the subgroup
$G_{({\bf e})}$ of elements fixing pointwise $(-\infty ,a]$ 
for some $a \le a_{0}$. \parskip0pt

The proof of Lemma \ref{1.10} shows that the condition
$h_{1}\prec h_{2}$ means the existence of $a \in L_{g}$
with $h_{1}(a') < h_{2}(a')$ for all $a'<a$.
This obviously implies $\nu_{a_{0}}(h_{1})< \nu_{a_{0}}(h_{2})$.
We see that $(L^{a_{0}}_{g},<)$ can be identified with
the group $(G_{{\bf e}}/G_{({\bf e})},\prec )$. $\square$

\begin{lem} \label{1.12}
If the ordering of $L^{a_{0}}_{g}$ is not dense, then
$L^{a_{0}}_{g}$ is a cyclic group with respect to
the structure of $G_{{\bf e}}/G_{({\bf e})}$.
\end{lem}

{\em Proof.} Notice that if there is an interval $(a,b)$,
$a,b \in L^{a_{0}}_{g}$, which does not have elements from
$L^{a_{0}}_{g}$, then every $c\in L^{a_{0}}_{g}$ has
a successor from $L^{a_{0}}_{g}$.
This follows from the fact that $[a,g(a)]$ and $[c,g(c)]$ can
be taken onto $[g^{-1}(a_{0}),a_{0}]$ by an element from
$G_{{\bf e}}$.
$\square$

\section{Polish groups with comeagre conjugacy classes}

Dugald Macpherson and Simon Thomas have proved in \cite{mactho}
that if a Polish group has a comeagre conjugacy class then
every element of the group fixes a point under any action on
a $\mathbb{Z}$-tree without inversions.
In this section we generalize that result to the situation
which covers the case of isometric actions.
Our method is different and is based on some algebraic property
of comeagre classes (Condition (1) of Proposition \ref{4.1}).
As a result the theorem can be formulated in 
elementary terms not involving Polish groups.
We apply this to actions of the group $Sym(\omega)$.
In the second part of the section we study actions of 
groups with comeagre conjugacy classes and invariant ends.

\subsection{Fixed points}
In the case of an isometric action if $h_{2}\cdot h_{1} = h_{3}$
and $h_{1},h_{2},h_{3}$ are loxodromic, then Proposition \ref{3loxo}
implies that $h_1 ,h_2 ,h_3$ and $h^{-1}_2 h_1$ can not belong to 
the same conjugacy class.
Indeed, if for example $[b,h_{2}(b)]$ is a proper subsegment of
$[c,h_{3}(c)]$, where $b \in L_{h_{2}}$, $c \in L_{h_{3}}$, then
for no $c' \in L_{h_{3}}$ the segment $[c',h_{3}(c')]$ can be
mapped to $[b,h_{2}(b)]$ by a map induced by some $g \in G$.
\parskip0pt

This observation motivates the following proposition.

\begin{prop} \label{4.1}
Let a group $G$ have a classical action on a pretree $T$.
Let $X \subset G$ satisfy the following conditions.\parskip0pt

(1) For every sequence $g_{1},...,g_{m} \in G$ there exist
$h_{0},h_{1},...,h_{m} \in X$ such that for every $1 \le i \le m$,
$g_{i} = h_{0}h_{i}$. \parskip0pt

(2) If $T^{h}_{0} = \emptyset$ for some $h \in X$, then all
$h \in X$ are loxodromic and there are no $h_{1},h_{2}\in X$
and $c_{1} \in L_{h_{1}}, c_{2} \in L_{h_{2}},$
such that $[c_{1},h_{1}(c_{1})]$ properly contains
$[c_{2},h_{2}(c_{2})]$. 
\parskip0pt

Then for any $g\in G$, $T^{g}_{0}\not= \emptyset$.
\end{prop}

{\em Proof.}
If all $h \in X$ are loxodromic, find $h_{1},h_{2},h_{3},h_4 \in X$
such that $h_{1} \cdot h_{3} = h_{2}$ and $h_{1} h_4 = h^{-1}_2$ (by (1)).
Then $h^{-1}_4 =h_2 h_1$ , $h^{-1}_3 =h^{-1}_2 h_1$ and by 
Proposition \ref{3loxo} there are $i,j \in \{ 1,2,3,4\}$, $i\not= j$, 
and $c \in L_{h_{i}}$ and $d \in L_{h_{j}}$ such that the segment
$[c,h_{i}(c)]$ properly contains $[d,h_{j}(d)]$.
This refutes condition (2). \parskip0pt

We now have that for all $h \in X$, $T^{h}_{0} \not= \emptyset$.
Let $g \in G$ have no fixed points in $T_{0}$ and
$P^{g} = \{ p_{1},p_{2}\}$ ($g$ is loxodromic by $(C2)$ of
Definition \ref{defcl}).
Let $g = h'\cdot h$, $h,h' \in X$.
Find $a$ and $b$ satisfying the properties that
$\{ a\} \cup T^{h}_{0}$ and $\{ b\} \cup T^{h'}_{0}$ are full
and for all $a' \in T^{h}_{0}$ and $b' \in T^{h'}_{0}$,
$\{ a,b\}\subset [a',b']$ (Lemma \ref{4sts}(3); note that formally 
it may happen that $a\not\in T^{h}_{0}$ or $b\not\in T^{h'}_{0}$).
By (C4) of Definition \ref{defcl} there are
$a_{0}\in T^{h}_{0}\cap [p_{1},p_{2}]$ and
$b_{0}\in T^{h'}_{0}\cap [p_{1},p_{2}]$.
This implies that $a, b \in [p_{1},p_{2}]$.\parskip0pt

By (1) there exist  $h_{0},h_{1},h_{2},h_{3} \in X$ such that
$g = h_{0}h_{3}$, $h = h_{0}h_{1}$, $h' = h_{0}h_{2}$.
By (C3) of Definition \ref{defcl} there are
$a_{1} \in T^{h_{0}}_{0} \cap T^{h}_{0}$
and $b_{1} \in T^{h_{0}}_{0} \cap T^{h'}_{0}$.
Then by (C0), $a,b \in  T^{h_{0}}_{0}$.
Applying (C0) again we see
$m(a,a_{0},a_{1}) \in T^{h_{0}}_{0} \cap T^{h}_{0}\cap [p_{1},p_{2}]$
and
$m(b,b_{0},b_{1}) \in T^{h_{0}}_{0} \cap T^{h'}_{0}\cap [p_{1},p_{2}]$.
On the other hand, by (C4) applied to $h_{0}$ and $h_{3}$ the
intersection $T^{h_{0}} \cap [p_{1},p_{2}]$ is a singleton;
then $m(a,a_{0},a_{1}) = m(b,b_{0},b_{1})$ and
$T^{h}_{0}\cap T^{h'}_{0} \not= \emptyset$ contradicting the
assumption that $g$ does not have fixed points in $T_{0}$.
$\square$
\bigskip

It is worth noting that condition (1) of Proposition 3.1 can be
weakened assuming that $m = 3$.\parskip0pt

A Polish group is a topological group whose topology is
{\em Polish} (a Polish space is a separable completely
metrizable topological space).
A subset is {\em comeagre} if it containes
an intersection of a countable family of dense open sets.

\begin{prop} \label{4.2}
Let $G$ be a Polish group with a classical action
on a pretree $T = T_{0} \cup P$.
Let $X \subseteq G$ be comeagre.
\parskip0pt

(1)  Then $X$ satisfies condition (1) of Proposition \ref{4.1}.
If $X$ is a conjugacy class then $X = X^{-1}$. \parskip0pt

(2) Let $X$ consist of elements fixing points in $T_{0}$.
Then every $g \in G$ has a fixed point in $T_{0}$.
\end{prop}

{\em Proof.}
(1) Let $g_{1},...,g_{m} \in G$.
Since the set $X$ is comeagre in  $G$, all $g_{i}X^{-1}$ are
comeagre and have a common element $h_{0}$.
Now find $h_{1},...,h_{m} \in X$ such that for
any $1 \le i \le m$, $g_{i} = h_{0}h_{i}$.
If $X$ is a conjugacy class then $X = g^{G}$, where
$g \in X \cap X^{-1}$.
Thus $X = X^{-1}$. \parskip0pt

(2) By (1) we can apply the proof of Proposition \ref{4.1}.
$\square$
\bigskip

We now see that Proposition \ref{4.1} generalizes the result of
Macpherson and Thomas mentioned above.
Indeed, let a Polish group $G$ have a comeagre conjugacy class
$X$.
By Proposition \ref{4.2} the group $G$ satisfies condition (1)
of Proposition \ref{4.1}.
If $G$ has an isometric action on an $\mathbb{R}$-tree then
condition (2) of Proposition \ref{4.1} is obvious.
Thus every element of $G$ fixes a point.

\subsection{Permutation groups}

Here we give an application of Proposition \ref{4.2}.
Let $A(\mathbb{Q})$ be the group of order-preserving
permutations of the rationals.
Then $A(\mathbb{Q})$ can be considered as a subgroup of
$Sym(\omega)$.
The following theorem shows that classical actions of
the symmetric group are determined by $A(\mathbb{Q})$.\parskip0pt

\begin{thm}
Let $Sym(\mathbb{Q})$ have a classical action on a pretree
$T = T_{0} \cup P$.
If for the corresponding action of $A(\mathbb{Q})$ every element
fixes a point from $T_{0}$, then so does every element
of $Sym(\mathbb{Q})$.
\end{thm}

{\em Proof.}
We define an expansion of the structure $(\mathbb{Q},<)$ by
relations
\footnote{denoting cycles of length $n$}
$(P_{n} : n \in \omega \setminus \{ 0\})$.
The expansion satisfies the  following properties:
$$
\forall x_{1},...,x_{n} (P_{n}(x_{1},...,x_{n}) \rightarrow
x_{1} <... < x_{n});
$$
$$
\forall x_{1},...,x_{n},y_{1},...,y_{m}
((\{ x_{1},...,x_{n}\} \cap \{ y_{1},...,y_{m}\} \not= \emptyset)  \wedge
P_{n}(x_{1},...,x_{n}) \wedge P_{m}(y_{1},...,y_{m}) \rightarrow
$$
$$
\{ x_{1},...,x_{n}\} = \{ y_{1},...,y_{m}\});
$$
$$
\forall n,y_{1},y_{2}\exists x_{1},...,x_{n}(y_{1}<y_{2}\rightarrow
P_{n}(x_{1},...,x_{n}) \wedge y_{1} < x_{1}<...< x_{n}<y_{2});
$$
$$
\forall x (\exists n,x_{1},...,x_{n})(P_{n}(x_{1},...,x_{n})\wedge
(x \in \{ x_{1},...,x_{n} \})).
$$
Such an expansion can be easily obtained using the fact that
the rationals form a countable dense linear ordering without
ends (then having an expansion where all conditions till the
last one are satisfied, put the elements for which the last
condition does not hold, into $P_{1}$).\parskip0pt

By back-and-forth
\footnote{see \cite{CK}: Ex.1.3.15 for the definition and
Theorem 1.4.2 for an illustration}
we now find an increasing $f \in A(\mathbb{Q})$ such that :
\parskip0pt

- each orbit of $f$ is cofinal in $\mathbb{Q}$; \parskip0pt

- the union of all $f$-orbits included in $P_{1}$ is dense in
$\mathbb{Q}$ and \parskip0pt

$$
\forall n, x_{1},...,x_{n} (P_{n}(x_{1},...,x_{n}) \rightarrow
(\exists k\in\mathbb{Z})(\mbox{ assuming that } B=\mathbb{Z}\setminus A
$$
$$
\mbox{ and } A =\{ i : k< i\} \mbox{ or } A =\{ i : i \le k \}
\mbox{ we have }
$$
$$
((\forall i \in A)P_{n}(f^{i}(x_{1}),...,f^{i}(x_{n})) \vee
(\forall i \in A)\bigwedge_{j\le n} P_{1}(f^{i}(x_{j}))) \wedge
$$
$$
((\forall i \in B) P_{n}(f^{i}(x_{1}),...,f^{i}(x_{n}))\vee
(\forall i \in B)\bigwedge_{j\le n} P_{1}(f^{i}(x_{j}))));
$$
$$
\forall n,y_{1},y_{2}\exists x_{1},...,x_{n},z_{1},...,z_{n}
(y_{1}<y_{2}\rightarrow (y_{1}<x_{1}<...<x_{n}<z_{1}<...<z_{n}
< y_{2} \wedge
$$
$$
(\forall i \le 0) (P_{n}(f^{i}(x_{1}),...,f^{i}(x_{n}))\wedge
\bigwedge_{j\le n} P_{1}(f^{i}(z_{j})))\wedge
$$
$$
(\forall i >0)(P_{n}(f^{i}(z_{1}),...,f^{i}(z_{n}))\wedge \bigwedge_{j\le n}
P_{1}(f^{i}(x_{j}))))).
$$
Let $h$ be the permutation of $\mathbb{Q}$ defined by:
$P_{n}(x_{1},...,x_{n})\rightarrow h(x_{1})=x_{2} \wedge ...\wedge
h(x_{n}) = x_{1}$.
It is easily seen that the permutation $g = fh^{-1}f^{-1}h$
has infinitely many cycles of each length.  \parskip0pt

Define a permutation $h'$ as follows.
For each $n > 1$ and $(x_{1},...,x_{n}) \in P_{n}$ with
$\bigwedge P_{1}(f(x_{j}))$ replace every $h$-cycle
$f^{-i}(x_{j}),j \le  n,$ with odd $i>0$ by $n$ single cycles
and create $n$-element cycles
$f^{i}(x_{1})\rightarrow ...\rightarrow f^{i}(x_{n})\rightarrow
f^{i}(x_{1})$ for all $i > 0$. \parskip0pt

If $\bigwedge P_{1}(f^{-1}(x_{j})),j \le n,$ holds for
$(x_{1},...,x_{n}) \in P_{n}$ (by the definition of $f$ this is
incompatible with the situation of the previous paragraph),
then create $n$-element cycles
$f^{-i}(x_{1})\rightarrow ...\rightarrow f^{-i}(x_{n})\rightarrow
f^{-i}(x_{1})$ for even $i> 0$ and remove all  $h$-cycles
$f^{i}(x_{j}),j \le  n,$ for $i \ge 0$.
We now see that $h'f^{-2}(h')^{-1}(x) = f^{-1}h^{-1}f^{-1}h(x)$.
\parskip0pt

As a result the permutation $g$ is the product of $f$
and $h^{-1}f^{-1}h$, where $f \in A(\mathbb{Q})$, and 
the permutation $h^{-1}f^{-1}h$ is  the product of $f$ and
$h'f^{-2}(h')^{-1}$ ($ = f^{-1}h^{-1}f^{-1}h$). \parskip0pt

If $Sym(\mathbb{Q})$ has a classical action on a pretree $T$,
then by the assumption each element of a conjugacy class meeting
$A(\mathbb{Q})$ fixes a point of $T_{0}$.
Now by $(C3)$, $f$ and $h^{-1}f^{-1}h$ have a common fixed point
in $T_{0}$, which is a fixed point of $g$, so $g$ is not
loxodromic. \parskip0pt

The permutation $g$ represents the comeagre conjugacy class in
$Sym(\omega)$, \cite{truss2}.
Now  Proposition \ref{4.2} (assuming that $X$ is the comeagre
conjugacy class) works in our case. $\square$

\subsection{Comeagre conjugacy classes and invariant ends}

In the following proposition we consider a situation
which appears in the case when a Polish group acts
with an invariant end.

\begin{prop} \label{4.4}
Let a group $G$ have a non-nesting action on an
$\mathbb{R}$-tree $T_{0}$ with an invariant end.
Let $X \subset G$ be a conjugacy class with $X^{-1} = X$ and
the following condition:
\begin{quote}
For  any $g_{1},g_{2},g_{3} \in G$ there exist
$h_{0},h_{1},h_{2},h_{3} \in X$ such that for any $1\le i\le 3$,
$g_{i} = h_{0}h_{i}$.
\end{quote}
Then for any $g \in G$, $T^{g}_{0} \not= \emptyset$.
\end{prop}

{\em Proof.}
If $T^{h}_{0} \not= \emptyset$ for some $h \in X$,
then all $h \in X$ fix some points.
This case can be considered as in the proof of
Proposition \ref{4.1}.  \parskip0pt

If $T^{h}_{0} = \emptyset$ for some $h \in X$, then all
$h \in X$ are loxodromic.
We want to show that this case is impossible.
Let a half-line $[t_{0},\infty)$ represent the invariant end.
\parskip0pt

Find $h_{0},h_{1},h_{2}\in X$ with  $h_{0}\cdot h_{1} = h_{2}$
and $a\in [t_{0},\infty)$ belonging to the axes
of these elements.
Replacing appropriate $h_{i}$ by $h^{-1}_{i}$ if necessary
(and moving elements from one side to another),
we may assume that all $h_{i}$ are increasing on $[a,\infty)$.
Let $h_{0}(a) \le h_{1}(a)$ (the case $h_{0}(a) > h_{1}(a)$ is
similar) and let $g \in G$ satisfy $gh_{0}g^{-1} = h_{2}$.
Since $g$ fixes the same end with $h_{0},h_{1},h_{2}$, it must
be loxodromic (otherwise $h_{0}$ and $h_{2}$ eventually coincide
on $[t_{0},\infty )$ and $h_{1}$ is not loxodromic).
\parskip0pt

Assume that $L_{g}$ contains $[a,\infty)$.
Find $b_{1},b_{2} \in [a,\infty)$ such that $g(b_{1}) = b_{2}$.
Assume $b_{1} \le b_{2}$.
Let $m$ be the minimal number such that
$h^{m}_{0}(b_{1}) > b_{2}$ (notice $1 \le m$).
Non-nesting and the condition $h_{0}(a) \le h_{1}(a)$ imply
$h^{m}_{0}(b_{1}) \le h_{1}(b_{2})$.
Then  $h^{m}_{0}g^{-1}(b_{2}) = h^{m}_{0}(b_{1}) > b_{2}$ and
$$
h^{m}_{0}g^{-1}(h_{2}(b_{2})) = h^{m+1}_{0}(b_{1}) =
h_{0}h_{0}(h^{m-1}_{0}(b_{1})) \le h_{0}h_{1}(b_{2}) = h_{2}(b_{2}).
$$
We now see that the element $h^{m}_{0}g^{-1}$ maps
$[b_{2},h_{2}(b_{2})]$ properly into itself.
This is a contradiction. \parskip0pt

If $b_{2} < b_{1}$ let $m$ be the minimal number such that
$h^{-m}_{0}(b_{1}) \le b_{2}$ (notice $1 \le m$).
Non-nesting and the condition $h_{0}(a) \le h_{1}(a)$ imply
$h^{-m+1}_{0}(b_{1}) \le h_{1}(b_{2})$.
Then $h^{-m+1}_{0}g^{-1}(b_{2})= h^{-m+1}_{0}(b_{1})> b_{2}$ and
$$
h^{-m+1}_{0}g^{-1}(h_{2}(b_{2})) = h^{-m+2}_{0}(b_{1}) =
h_{0}(h^{-m+1}_{0}(b_{1})) \le h_{0}h_{1}(b_{2}) = h_{2}(b_{2}).
$$
We now see that the element $h^{-m+1}_{0}g^{-1}$ maps
$[b_{2},h_{2}(b_{2})]$ properly into itself.
This is a contradiction.
$\square$
\bigskip

We now conclude the material above by the following theorem.

\begin{thm} \label{TheoremA}
Let a Polish group $G$ have a non-nesting action on an
$\mathbb{R}$-tree $T_{0}$ without $G$-fixed points in $T_{0}$.
Let $X \subseteq G$ be a comeagre set.
Then the following statements hold. \parskip0pt

If every element of $X$ fixes a point, then every element of $G$
fixes a point. \parskip0pt

If $G$ fixes an end and $X$ is a conjugacy class, then
every element of $G$ fixes a point.
\end{thm}

{\em Proof.}
We already know that the assumptions imply that 
the action is classical (Proposition \ref{3.1}).
Now the first claim of the theorem follows from
Proposition \ref{4.2}(2).
By Proposition \ref{4.2}(1) and Proposition \ref{4.4} 
we have that if a Polish group $G$ has a comeagre conjugacy 
class then every element of $G$ fixes a point under any 
non-nesting action on an $\mathbb{R}$-tree with an invariant end.
This is the second claim of the theorem. $\square$
\bigskip

It is known that $Sym(\omega)$ and $A(\mathbb{Q})$ 
have comeagre conjugacy classes \cite{truss2}.

\section{Comeagre conjugacy classes and end stabilizers}

To formulate the main result of the section we need
the following definition.
For a subset $A$ of a median pretree $T$ define the 
{\em closure} $c(A)$ as the minimal subpretree of $T$ 
with the property that the $T$-median of any triple from 
$c(A)$ belongs to $c(A)$.
It is clear that in the case when a group $G$ acts on 
a pretree $T$ and $A$ is $G$-invariant, the pretree 
$c(A)$ is $G$-invariant too.  
\parskip0pt

The following theorem roughly says that the presence of
a comeagre loxodromic conjugacy class implies that the
$G_{{\bf e}}$-orbits are much smaller than the corresponding
$G$-orbits.

\begin{thm} \label{comp}
Let a group $G$ have a non-nesting action on
an $\mathbb{R}$-tree $T_{0}$.
Let $X \subset G$ be a conjugacy class of loxodromic
elements satisfying the following condition.\parskip0pt
\begin{quote}
For  every triple $g_{1},g_{2},g_{3} \in G$ there 
exist $h_{0},h_{1},h_{2},h_{3}\in X$ such that 
$g_{i} = h_{0}h_{i}$ for all $1\le i\le 3$.
\end{quote}
Then for any $g \in X$, an end ${\bf e}$ represented by
$L_{g}$ and a point $a_{0} \in L_{g}$ the ordering
$L^{a_{0}}_{g}=G_{{\bf e}}a_{0}\cap L_{g}$
is not dense in $c(Ga_{0})\cap L_{g}$.
In particular this conclusion holds if $X$ is
a comeagre conjugacy class of loxodromic elements.
\end{thm}

To illustrate some aspects of the formulation let $g\in G$ 
be loxodromic and $a_{0}\in L_{g}$ be as in the theorem.
It is clear that the set $L^{a_{0}}_{g}=G_{{\bf e}}a_{0}\cap L_{g}$ 
is cofinal (in both directions) in the line $L_{g}$.
On the other hand it may happen that the ordering
$c(Ga_{0})\cap L_g$ (induced by a natural ordering of $L_{g}$)
is dense but not dense in $L_{g}$.

\bigskip

{\bf Example.} Consider $\mathbb{R}$ as
$\mathbb{Z}\times \{ a,b\}\times (0,1]$, where the elements of
$(2k,2k+1]$ are denoted by triples $(k,a,r)$, $r \in (0,1]$ and
the elements of $(2k+1,2k+2]$ are denoted by triples $(k,b,r)$,
$r \in (0,1]$, $k\in \mathbb{Z}$.
Let $G = \mathbb{Q}$  act on $\mathbb{R}$ as follows.
If $q + k + r = k' + r'$ , with $k'\in \mathbb{Z}$ and
$r' \in (0,1]$, then put $(k,a,r)+q = (k',a,r')$ and
$(k,b,r)+q = (k',b,r')$.
It is easy to see that the action of $\mathbb{Q}$ obtained on
$\mathbb{R}$ is non-nesting.
On the other hand the interval $(0,1]$ (consisting of all
$(0,a,r)$ with $r \in (0,1]$) does not contain any element of
the orbit of $0= (-1,b,1)$.
The fact that the orbit of $0$ is a dense ordering follows from
density of $\mathbb{Q}$.
Its closure coincides with the orbit. $\square$
\bigskip

We now describe one of our tools.
A binary relation $r$ (a partial ordering, where
$r(a,b)\vee (a=b)$ is interpreted as $a\le b$) on a pretree 
$T$ is called a {\em flow} (\cite{bow}, pp. 23 - 25) if
it satisfies the following axioms: \parskip0pt

$\neg (r(x,y)\wedge r(y,x))$,
$B(z;x,y)\rightarrow r(x,z)\vee r(y,z)$ and \parskip0pt

$(r(x,y)\wedge z \not=y) \rightarrow (B(y;x,z)\vee r(z,y))$.

The material of the next paragraph is based on pp. 26 - 28 of 
\cite{bow}. 
It would be helpful for the reader (but not necessary) 
to recall some formulations given there. \parskip0pt

We say that a flow $r$ is
{\em induced by an endless directed arc} $(C,<)$ if
$(x,y)\in r\leftrightarrow\exists z\in C\forall w >z B(y;x,w)$
\footnote{by Lemma 3.8 of \cite{bow} for any endless directed arc $(C,<)$ 
this formula defines a flow}
(see \cite{bow}, p. 26).
Then it is easy to see that for any arc $J$ if $J$ does not have
maximal elements with respect to $r$, then the formula
$r(x,y)\vee r(y,x)\vee x = y$ defines an equivalence relation
on $J$ with at most two classes.
We say that $r$ {\em lies} on $J$ if $J$ contains a maximal
element of $r$ in $T$ or $J$ does not have $J$-maximal elements
with respect to $r$ and the equivalence relation
$r(x,y)\vee r(y,x)\vee x = y$ defines a non-trivial cut
$J= J^{-}\cup J^{+}$, such that for any $a \in J^{+}$,
$b\in J^{-}$ there is no $c \in T$ with $r(a,c)\wedge r(b,c)$
(so it may happen that $C\cap J$ is cofinal with $C$).
It is easy to see that for any line $L$ there is a natural 
function from the set of all flows of $T$ induced by endless 
directed arcs and lying on $L$ onto the set of all Dedekind 
cuts on $L$: the Dedekind cut corresponding to a flow $r$ is 
determined by a maximal element or the equivalence relation 
$r(x,y)\vee r(y,x)\vee x = y$.
The following lemma shows that when $T$ is median and dense, 
this correspondence is bijective for flows without maximal elements.

\begin{lem} \label{flowcut} 
Assume that $T$ is a median pretree and $r$ is a flow induced by 
an endless directed arc. 
Let $C,D$ be arcs of $T_0$ which are linearly ordered with respect 
to $r$ and do not have upper $r$-bounds in $T_0$.
Then $C$ and $D$ are cofinal or for any $c \in C$ and $d \in D$,
$[c,d]=\{ a: a\in C\wedge r(c,a)\vee a\in D\wedge r(d,a)\}$.
Each of these arcs defines the flow $r$ as in the definition above.
\end{lem}

{\em Proof.}
To see the first statement we start with the case when 
$r$ is induced by $C$.
Let $c$ and $d$ be as in the formulation. 
If $C$ and $D$ are not cofinal then $C\cap D=\emptyset$
(apply the fact that if $t\in C\cap D$ and $D\models t<t'$,
then $t'$ must belong to $C$).
By the definition of $r$, if $D\models d<d'$, then $d'$
belongs to some $[d,c']$ with $C\models c<c'$.
If $d'\not\in [d,c]$, then $d'$ belongs to $[c^{*},c']$, where
$c^{*}$ is the median of $c,c',d$, a contradiction with
$C\cap D=\emptyset$.
It is now easy to see that the set $[c,d]\cap D$ consists 
of all $d'\in D$ with $r(d,d')$. \parskip0pt

This implies that there is no $a\in [c,d]\setminus (C\cup D)$
(otherwise $a$ is an upper $r$-bound for $D$) and
there is no $c'\in C$ with $c'\not\in [c,d]$ and $c<c'$
(otherwise $c'$ is an upper $r$-bound for $D$).
Now a straightforward argument gives the formula for $[c,d]$ 
as above. \parskip0pt

In the case when $r$ is induced by some ordering $A$ and
$A$ is cofinal with $C$ or $D$, the argument above works again.
If $A$ is not cofinal with these orderings then for every $a\in A$
we have
$[c,a]=\{ a': a'\in C\wedge r(c,a')\vee a'\in A\wedge r(a,a')\}$
and
$[a,d]=\{ a': a'\in D\wedge r(d,a')\vee a'\in A\wedge r(a,a')\}$.
Then the median of $a,c,d$ must belong to one of the intervals $C$, 
$D$ or $A$. 
The case $m(a,c,d)\in A$ is impossible, because the elements of $A$ 
greater than $m(a,c,d)$ cannot belong both to $[c,a]$ and $[d,a]$. 
When $m(a,c,d)\in C\cup D$, the arcs $C$ and $D$ are cofinal. 
\parskip0pt

To show that {\em the flow $r$ is induced by
any of its linear orderings without upper $r$-bounds in $T_0$}
we apply similar arguments as above.
Indeed, if $A$ and $C$ are linear orderings without upper $r$-bounds
in $T_0$, $r$ is induced by $A$ and $C$ is not cofinal with $A$,
then any $D$ as above is cofinal with $A$ or $C$ 
(as $T$ is median).
If $D$ is cofinal with $A$, then $D$ induces $r$.
If $D$ is cofinal with $C$, then for any $d\in D$ and $a\in A$,
$[a,d]=\{ a':a'\in D\wedge r(d,a')\vee a'\in A\wedge r(a,a')\}$.
Let $t\in T_0$ and $t^{*}=m(t,a,d)$ for some $d\in D$ and $a\in A$.
Now it is straightforward that for any $t'\in T_0$ the condition
$r(t,t')$ is equivalent to
$t'\in (t,t^{*}]\vee t'\in (t^{*},d]\cap A\vee t'\in (t^{*},a]\cap D$.
Using this formula it is easy to verify that $r$ is induced by $D$ 
(notice that $A$ and $D$ are symmetric in this condition). $\square$ 
\bigskip 

By Lemma \ref{flowcut} we see that when $r$ is a flow of 
a dense median pretree $T$ defined by an endless directed arc  
and $r$ lies on a line $L$ but does not have a maximal element 
lying on $L$, then each of the half-lines on $L$ defined by the 
corresponding equivalence relation is an endless directed arc 
inducing $r$.  

\bigskip

{\em Proof of Theorem \ref{comp}.}
By Theorem \ref{TheoremA} we may assume that $T_0$ is not a line.
If the theorem is not true there are $g \in X$,
a point $a_{0}\in L_{g}$ and an end ${\bf e}$
represented by $L_{g}$ such that the ordering
$L^{a_{0}}_{g} = G_{{\bf e}}a_{0}\cap L_{g}$ is dense in
the line $c(Ga_{0})\cap L_{g}$ of the pretree $c(Ga_{0})$.
Note that $c(Ga_0 )\cap L_g$ is cofinal in $L_g$ and this 
implies the same statement for any line $L_h$ with $h\in X$.
\parskip0pt

It is also worth noting that the subspace
$\bigcup\{ L_{h}:h\in g^{G}\}$ is a full subtree of $T_{0}$.
Indeed, for any $h,h'\in X$ there exist $h_{0},h_{1},h_{2}\in X$
such that $h = h_{0}h_{1}$ and $h' = h_{0}h_{2}$.
By Proposition \ref{3loxo} there is an arc in
$L_{h_{0}}\cup L_{h_{1}}\cup L_{h_{2}}$ joining $L_{h}$
and $L_{h'}$.
{\em We may assume that} $T_{0}=\bigcup\{ L_{h}:h\in X\}$.
\parskip0pt

Notice that if $a\in L_{h'} \cap L_{h''}$ 
defines an $L_{h'}$-half-line without other common elements 
with $L_{h''}$, then $a$ is the median of three non-linear 
elements from $Ga_0$ and thus belongs to $c(Ga_0 )$.
In particular in the situation above the arc joining $L_h$ 
and $L_{h'}$ consists of at most three intervals with 
extremities from $c(Ga_0 )$.   \parskip0pt 
   
{\em We want to embed the $G$-pretree $c(Ga_0 )$ into some 
special $\mathbb{R}$-tree with an isometric action of $G$.} 
We start with the case when $L^{a_0}_g$ is not a dense ordering. 
Here we apply the results of Section 2.2 and consider 
$L^{a_0}_g$ as the ordered group $G_{{\bf e}}/G_{({\bf e})}$. 
In this case all $h^{-1}(L^{a_0}_g )$ are discrete. 
Since $T_0 =\bigcup \{ L_h :h\in X\}$ we see that $c(Ga_0 )$ 
is discrete. 
Thus $c(Ga_0 )$ can be considered as a simplicial tree with 
an isometric action of $G$. 
Now the pretree $c(Ga_0 )$ can be naturally 
\footnote{as any simplicial tree}
embedded into an $\mathbb{R}$-tree with an isometric action of $G$. 
We have obtained a classical action of $G$ on a tree satisfying 
the conditions of Proposition \ref{4.1} (the second one holds 
because the action is isometric).
By Proposition \ref{4.1} the elements of $X$ are not loxodromic, 
a contradiction.\parskip0pt   

From now on we consider the case when $L^{a_0}_g$ is dense. 
Denote $c(Ga_0 )$ by $T'$.
If elements $a_{1}$ and $a_{2}$ belong to $T'$, say 
$a_{1}\in L_{h}$ and $a_2 \in L_{h'}$, 
then as above we find $h_{0},h_{1},h_{2}\in X$ such that 
there is a $T_0$-arc in $L_{h_{0}}\cup L_{h_{1}}\cup L_{h_{2}}$ 
joining $L_{h}$ and $L_{h'}$.
Since $c(Ga_0 )$ is median, the intervals of the corresponding 
lines have extremities belonging to $c(Ga_0 )$.
This implies that the $T'$-interval $[a_1 ,a_2 ]$ is decomposed 
in $T'$ into at most five intervals from the corresponding lines. 
We now see that all intervals of $T'$ are dense.
Let $T^{*}$ be the set consisting of $T'$ and all flows
of $T'$ which are induced by endless directed arcs and 
which do not have maximal elements.
We will show below that $T^{*}$ can be presented as an 
$\mathbb{R}$-tree where the action of $G$ is isometric.  
We start with some helpful observation. \parskip0pt 

{\bf Claim 1.} {\em Every endless directed arc $I$ from $T'$
which does not define an end of $T_0$, is cofinal with
an endless directed arc from some $L_{h}\cap T'$, $h\in X$.}
\parskip0pt 

Indeed, let $a\in I$.
Since $T_{0}$ together with the set of ends forms a complete 
tree, there is $c\in T_{0}$ such 
that $I$ is cofinal with $[a,c)\cap c(Ga_0 )$.
If $a\in L_{h'}$ and $c\in L_{h}$, then as above we find
$h_{0},h_{1},h_{2}\in X$ such that there is an arc in
$L_{h_{0}}\cup L_{h_{1}}\cup L_{h_{2}}$ joining $L_{h}$
and $L_{h'}$.
We see that $[a,c)$ is decomposed into at most five
intervals from the corresponding lines.
The last interval (which is of the form $[a',c)$ with 
$a'\in c(Ga_0 )$) is cofinal with $I$. \parskip0pt

We now extend the betweenness relation to $T^{*}$ as in \cite{bow}
(p.30): for points $x,y$ and a flow $r$ we say $B(y;x,r)$ if
$(x,y) \in r$.
Then we can define $B(r;x,y)$ as the case when
there is no point $z$ with $B(z;x,r)\wedge B(z;y,r)$.
For flows $r,r'$ and a point $x$ we say $B(x;r,r')$ if for
any point $y\not= x$, $(y,x)\in r\vee (y,x)\in r'$.
We also put $B(r;x,r')$ if $B(r;x,y)$ for some $y\in T'$ with
$(x,y)\in r'$.
If $p,q,r$ are flows then $B(p;q,r)$ means that there is no
point $z$ with $B(z;p,q)\wedge B(z;p,r)$. \parskip0pt 

This definition implies that for {\em a non-linear triple
$u,v,w\in T^{*}$ the median $m(u,v,w)$ can be presented
as $m(x,y,z)$ for some non-linear $x,y,z\in T'$;
thus $m(u,v,w)\in T'$}.\parskip0pt 

By Claim 1 we know that {\em all flows of $T^* \setminus T'$ 
are induced by endless directed arcs of lines} $L_h \cap T'$, 
$h\in X$. \parskip0pt 

It is also worth noting that the definition of the betweenness 
relation on $T^*$ implies that a flow $r\in T^* \setminus T'$ 
{\em lies on a line $L'\subset T'$ if and only if $r$ forms 
a linear triple with any two points of} $L'$.\parskip0pt 
 
As a flow $r\in T^{*}\setminus T'$ is induced by any of its
linear orderings without upper bounds (Lemma \ref{flowcut}), 
if $r$ belongs to a $T^{*}$-line $L^{*}$, then $r$ is defined 
by a class of the corresponding Dedekind cut of $L'=T'\cap L^{*}$.
Moreover any cut $L'= C^{-}\cup C^{+}$, which does not
define an element of $T'$, defines a flow from $T^{*}$.
Indeed suppose for a contradiction that $r$ is the flow 
defined by $C^{-}$ and $c$ is a maximal element of $r$.
Then for any $c'\in C^{-}$ and $c''\in C^{+}$,
the interval $(c',c)$ (in $T'$) consists of elements 
of $L'$ $r$-greater than $c'$ (by axioms of flows and 
maximality of $c$) and is contained in $[c',c'']$ 
(see Lemma \ref{flowcut}). 
Since $T'$ is quasimedian, $c\in [c',c'']$,
contradicting the assumption that the cut does not 
define an element of the line.
As a result we obtain that {\em the $T^{*}$-line
$L^{*}$ is the Dedekind completion of $L'$ and can be
identified with} $\mathbb{R}$. \parskip0pt

Since $T'$ is dense, $T'$ {\em is dense in} $T^*$.\parskip0pt 
 
The action of $G$ on $T'$ uniquely defines an action on $T^{*}$.
We want to show that this action is a non-nesting action on
an $\mathbb{R}$-tree.
Let us start with the following claim. \parskip0pt 

{\bf Claim 2.} {\em If $g'\in G$ is not loxodromic in $T_0$, 
then $T'$ contains a point fixed by $g'$ or a segment 
which is inversed by $g'$.} \parskip0pt 

Indeed, if $g'$ has a fixed point $a\in T_0$ such that 
for some $h_1 ,h_2 ,h_3 \in X$ the point $a$ defines 
three pairwise disjoint half-lines 
(without the extremity, half-lines of the form $(a,\infty )$) 
on the corresponding $L_{h_1}$, $L_{h_2}$ and $L_{h_3}$, 
then $a$ is the median of three elements of $Ga_0$ 
and thus belongs to $c(Ga_0 )$.
In particular if $g'$ fixes pointwise a segment or a half-line 
of some $L_h \not=g'(L_h )$, then $g'$ has a fixed point 
in $c(Ga_0 )$.    \parskip0pt 

Using non-nesting we see that in the remaining case we must 
consider the situation when $a$ with $g' (a)=a$ belongs to some 
line $L_h$, $h\in X$, and to some segment $[b,g' (b)]$ with
$b\in Ga_0$ and no point of $[b,a)$ is fixed by $g'$.
Since $a$ is not a median of three non-linear points in $T_0$,
we may assume that $b\in L_h$. 
By the same reason the interval $(a,g' g' (b))$ meets 
$(a,b)\cup (a,g' (b))$.
Our assumptions (together with non-nesting) imply that 
there is $c\in (a,b)\cap (a,g' g' (b))$ and this $c$ can be  
found in $c(Ga_0 )$.
By non-nesting, $g' g' (c)=c$. 
Thus $[c,g' (c)]$ is inversed by $g'$. \parskip0pt

{\bf Claim 3.} {\em An element $g'\in G$ fixes a point in $T_0$ 
if and only if it fixes a point in $T^{*}$.} \parskip0pt

If $g'$ does not fix a point in $T_0$, then there is a line 
$L\subseteq T_0$ which is the axis of $g'$.
Since every segment of $T_0$ can be presented as the union 
of at most five segments from lines $L_h$, $h\in X$, with 
common extremities belonging to $c(Ga_0 )$, 
the set $L\cap c(Ga_0 )$ is not empty and thus is cofinal with $L$. 
By the definition of $T^{*}$ and the corresponding betweenness 
relation, $L\cap Ga_0$ is cofinal with some $T^{*}$-line 
$L^{*}$.
It is clear, that $L^{*}$ is the axis of $g'$ in $T^{*}$.
Thus $g'$ is loxodromic in $T^{*}$ (by Lemma \ref{loxo}). 
\parskip0pt 

To see the converse we apply Claim 2.
By this claim we only have to consider the case when $g'$ 
fixes a point in $T_0$, does not fix any point of $T'$ and 
inverses a segment, say $[c,g'(c)]$, in $T'$. 
Let $I=\{ x\in [c,g'(c)]\cap T' :B(x;c,g'(x))\}$.
Since $T'$ is dense, $I$ is an endless directed arc. 
It is straightforward, that $I$ and $g'(I)$ define the same 
flow which is fixed by $g'$.     
\parskip0pt 

{\bf Claim 4.} {\em The action of the group $G$ on $T^{*}$ is 
non-nesting.}\parskip0pt 

Indeed, if $g'\in G$ is loxodromic in $T_0$, then it is loxodromic
in $T^{*}$. 
Lemma \ref{loxo} implies that if $g'$ maps an segment of $T^{*}$ 
properly into itself, then such a segment can be chosen in the 
axis $L^{*}_{g'}$ of $g'$. 
The latter is impossible, because $g'$ has a non-nesting action 
on $L^{*}_{g'} \cap T'$ with a cofinal orbit and $L^{*}_{g'}\cap T'$ 
is dense in $L^{*}$.   
If $g'$ fixes a point $a\in T^{*}$ and maps a segment $[b,b']$ 
properly into itself then $a$ and $[b,b']$ can be chosen so 
that $a\in [b,b']$. 
This can be shown by straightforward arguments depending 
on the place where the median of $a$ and the extremities of 
the segment lie.\parskip0pt 

If it happens that $a=b'$, then as $T'$ is dense in $[a,b]$, 
we can arrange that $b\in T'$.
Since the action of $G$ on $T'$ is non-nesting we see that 
$a\not\in T'$.
Using Claims 2 and 3 find a segment $[c,g'(c)]$ with 
$(g')^{2}(c)=c\in c(Ga_0 )$ and $a\in [c,g'(c)]$.
Replacing $c$ by $m(a,b,c)$ or $m(a,b,g'(c))$ if necessary,
we may assume that $c\in [a,b]$. 
Then $g'$ maps $[g'(c),b]$ properly into itself, contradicting 
non-nesting on $T'$. \parskip0pt 

Consider the case when neither $b$ nor $b'$ are fixed by $g'$.
If $g'(b)\in [a,b]$ or $g'(b')\in [a,b']$ then apply 
the argument of the previous paragraph.      
Assume $g'(b)\in [a,b']$. 
Then $g'(b')\in [a,b]$ and $(g')^2$ maps the segment $[a,b']$ 
properly into itself.
As we already know this contradicts the assumption that the 
action of $G$ on $T'$ is non-nesting. 
This finishes the proof of the claim.  
\parskip0pt

We now consider the action $*_{g}$ of $G_{{\bf e}}$ on 
$L_g\subseteq T_0$ (as in Section 2.2), where $g$ is as 
in the formulation of the theorem. 
Since every $g'\in G_{{\bf e}}$ acts on $L_g$ as a translation
(up to topological equivalence with $\mathbb{R}$) so it does  
on the corresponding axis $L^{*}_g$ from $T^{*}$
(by non-nesting). 
It is clear that all elements of $G_{({\bf e})}$ 
fix $L^{*}_g$ pointwise.     \parskip0pt 
          
As we already know, the line of $T^{*}$ containing a copy
$t^{-1}(L^{a_{0}}_{g})$ can be identified with
the Dedekind completion of the $T'$-line  
$t^{-1}(c(Ga_{0})\cap L_{g})$.
Since $L^{a_0}_{g}$ is dense in $c(Ga_0 )\cap L_g$, this 
completion coincides with the Dedekind completion of 
$t^{-1}(L^{a_0}_g )$. 
We will denote this line by $t^{-1}(L^{a_0}_{g})^{D}$.
It is clear that $t^{-1}(L^{a_0}_g )^{D}$ is the axis of 
$t^{-1}gt$ in $T^*$.
\parskip0pt

The line $L^{*}_g=(L^{a_{0}}_{g})^{D}$ can be identified 
with $\mathbb{R}$ so that the group $L^{a_{0}}_{g}$ 
(under the structure of the Archimedean group 
$G_{{\bf e}}/G_{({\bf e})}$ defined in Section 2.2) acts on
$(L^{a_{0}}_{g})^{D}$ by translations defined by real numbers.
We assume that $a_{0}$ corresponds to $0$.
Let $d$ denote the metric obtained on $(L^{a_{0}}_{g})^{D}$ in this way.
\parskip0pt

Now for any $h\in G$ consider $h^{-1}(L^{a_{0}}_{g})^{D}$ under
the metric $d^{h}$ induced by $d$ and the map $h^{-1}$ (so that
$h^{-1}$ is an isometry). \parskip0pt

{\bf Claim 5.} {\em The metrics $d^{t}$ and $d^{h}$ agree on
every common segment of the corresponding lines
\footnote{if the segment contains more than one element, then it 
contains infinitely many elements both fom $h^{-1}(L^{a_0}_g )$ 
and $t^{-1}(L^{a_0}_g )$}
.
In particular, if $(h^{-1}(L^{a_{0}}_g ))^D =(t^{-1}(L^{a_{0}}_g ))^D$, 
then $d^t =d^h$. } \parskip0pt

Suppose
$x_{1},x_{2}\in t^{-1}(L^{a_{0}}_{g})^{D}\cap h^{-1}(L^{a_{0}}_{g})^{D}$,
$x_{1} < x_{2}$ and $d^{t}(x_{1},x_{2}) < d^{h}(x_{1},x_{2})$
(the case $d^{h}(x_{1},x_{2}) < d^{t}(x_{1},x_{2})$ is similar).
Since $t^{-1}(L^{a_{0}}_{g})$ is dense in
$t^{-1}(L^{a_{0}}_{g})^{D}$, we can find
$x'_{1}\le x''_{1} < x''_{2}\le x'_{2}\in t^{-1}(L^{a_{0}}_{g})^{D}\cap h^{-1}(L^{a_{0}}_{g})^{D}$
with  $d^{t}(x'_{1},x'_{2}) < d^{h}(x''_{1},x''_{2})$
where $x'_{1},x'_{2} \in t^{-1}(L^{a_{0}}_{g})$ and
$x''_{1},x''_{2} \in h^{-1}(L^{a_{0}}_{g})$.
Then the inequality $t(x'_{2})-t(x'_{1})<h(x''_{2}) - h(x''_{1})$
holds in the group $L^{a_{0}}_{g}$.
Applying a translation from $L^{a_{0}}_{g}$ we can take
$h(x''_{1})$ to $t(x'_{1})$; then $h(x''_{2})$ must go to an
element of $L^{a_{0}}_{g}$ greater that
$t(x'_{2})$, a contradiction with non-nesting. \parskip0pt

Extend the metric $d$ to the space $T^*$ as follows.
Since $T'$ consists of
$\bigcup\{ h^{-1}(c(Ga_{0})\cap L_{g}): h\in G\}$ and
flows are defined by arcs from lines of the form
$h^{-1}(c(Ga_{0})\cap L_{g})$, any $a\in T^*$ belongs to 
some $T^*$-line of the form $h^{-1}(L^{a_{0}}_{g})^{D}$.
If $a\in h^{-1}(L^{a_{0}}_{g})^{D}$ and
$b\in (h')^{-1}(L^{a_{0}}_{g})^{D}$ then find
$h_{0},h_{1},h_{2} \in X$ such that $h = h_{0}h_{1}$ and
$h' = h_{0}h_{2}$.
There is a bridge
$[a',c_{1}]\cup [c_{1},c_{2}]\cup [c_{2},b']$ where
$a'\in h^{-1}(L^{a_{0}}_{g})^{D}$,
$b'\in (h')^{-1}(L^{a_{0}}_{g})^{D}$
and every segment above belongs to an appropriate line
$h^{-1}_{1}(L^{a_{0}}_{g})^{D}$,
$h^{-1}_{0}(L^{a_{0}}_{g})^{D}$ or $h^{-1}_{2}(L^{a_{0}}_{g})^{D}$
(it can happen that $a'=c_{1}$ or $c_{1}=c_{2}$ or $c_{2}= b'$).
Now define the distance between $a$ and $b$ as the sum of
distances in the sequence $a,a',c_{1},c_{2},b',b$ where each
distance is taken from the corresponding axis. \parskip0pt

By Claim 5 to prove that $d$ is invariant under the $G$-action 
it suffices to show that if $h$ maps $h^{-1}_{1}(L^{a_{0}}_{g})^{D}$
onto $h^{-1}_{2}(L^{a_{0}}_{g})^{D}$, then $h$ maps the metric
$d^{h_{1}}$ onto $d^{h_{2}}$.
The latter condition can be verified as follows. 
Since $h' =h_2 hh^{-1}_1$ preserves the line $L_g$, 
$h' (L^{a_0}_g )=L^{a_0}_g$.  
This shows that $h^{-1}$ maps $h^{-1}_2 (L^{a_0}_g )$ onto 
$h^{-1}_1 (L^{a_0}_g )$. 
By Claim 5 the metrics $d^{h_1}$ and $d^{h_2 h}$ coincide on 
$h^{-1}_1 (L^{a_0}_g )^D$.
The latter one is the image of $d^{h_2}$ under $h^{-1}$
(by the definition).  
\parskip0pt

The definition of $d$ implies that $T^*$ satisfies the definition 
of an $\mathbb{R}$-tree (see Section 2.1).
Moreover every line $h^{-1}(L^{a_{0}}_{g})^{D}$ becomes an
$\mathbb{R}$-line where $g^{h}$ acts by a translation. \parskip0pt

We have obtained a classical action on a tree satisfying
the conditions of Proposition \ref{4.1}.
By Proposition \ref{4.1} the elements of
$X$ are not loxodromic, a contradiction. $\square$

\end{document}